\documentclass[12pt,leqno]{amsart}
\usepackage{amsmath,amsthm}
\usepackage{amsfonts}

 \newtheorem{res}{Result}[section]
 \newtheorem{Ex}{Example}[section]
 \newtheorem{theorem}[res]{Theorem}
 \newtheorem{remark}[res]{Remark}
 \newtheorem{prop}[res]{Proposition}
 \newtheorem{lem}[res]{Lemma}
 \newtheorem{cor}[res]{Corollary}
 \newtheorem{defi}[res]{Definition}

\numberwithin{equation}{section}

\def\hatpi{\hat \Pi}
\def\defto{\buildrel def\over =}
\def\F{\mathcal F}
\def\calP{\mathcal P}
\def\K{\mathcal K}
\def\E{\mathbb E}
\def\A{\mathcal A}
\def\Q{\mathcal Q}
\def\Om{\Omega}
\def\om{\omega}
\def\al{\alpha}
\def\be{\beta}

\def\1{\mathbf 1}
\def\calF{\mathcal F}
\def\bbF{\mathbb F}
\def\R{\mathbb R}
\def\calI{\mathcal I}
\def\bbI{\mathbb I}
\def\P{\mathbb P}

\def\bbQ{\mathbb Q}
\def\vare{\varepsilon}
\def\to{\rightarrow}
\def\G{\mathcal G}
\def\Lam{\Lambda}
\def\undrho{\underline{\rho}}
\def\bbG{\mathbb G}
\def\Linf{L^{\infty}}
\def\Lplus{L^{\infty}_+}
\def\calH{\mathcal H}
\def\bbT{\mathbb T}

\begin{document}
\bibliographystyle{plain}
\title[On decomposing risk and reserving]{On decomposing risk\\
in a Financial-Intermediate market and reserving.}
\date{}
\author{SAUL JACKA}
\address{Dept. of Statistics, University of Warwick, Coventry CV4 7AL, UK}
\email{s.d.jacka@warwick.ac.uk}
\author{ABDELKAREM BERKAOUI}
\email{a-k.berkaoui@warwick.ac.uk}
\begin{abstract}
We consider the problem of decomposing monetary risk in the presence of
a fully traded market in {\it some} risks. We show that a mark-to-market
approach to pricing leads to such a decomposition if the risk measure
is time-consistent in the sense of Delbaen.
\end{abstract}



\thanks{This research was supported by the grant \lq Distributed Risk Management'
in the Quantitative Finance initiative funded by EPSRC and the
Institute and Faculty of Actuaries} \maketitle

\section{Introduction.}
In many contexts, financial products are priced and sold in the
absence of a market (i.e a fully traded market) in these products.
Typically these products have a dependence (either explicit or
implicit) on one or more securities or contracts in a traded
(financial) market. An obvious example is insurance (and, in
particular, life insurance), but other examples include (the
benefits provided by) pension funds and stock and options in
non-quoted companies. This paper is concerned with the questions of
valuing the liabilities of an (intermediate) market maker in such
products and of how to make and invest financial reserves for them.

We take the view that such a market maker is a price-taker in the
traded financial market (hereafter referred to simply as the market)
and is a price-maker in its own products (hereafter referred to as
contracts)---the ultimate value of which are contingent on risks not
present in the market. From this point of view we may add in other
non-market risk such as, for example, interruption of business,
fraud, litigation, insurable risks and economic factors (such as the
behavior of price indices and salaries) which impinge upon the
eventual settlement value (or payoff) of these contracts.

We adopt the view that in such a setup, the intermediate market
maker (hereafter referred to as the intermediate) will adopt a
coherent risk measure (on discounted final values) as their
valuation method and show how this implies certain constraints on
the form of this risk measure and finally, how if these constraints
are met, the risk measure implies a reserving method and an
investment strategy in the market.

\section{Contracts contingent on lives}
As we have already seen in the introduction, the issues we address
are by no means limited to life assurance and related products,
indeed they have relevance to monetary risk management in any
conceivable context; nevertheless, historically, life assurance and
annuities (the two main products of life insurance companies) are a
major source of such issues. Consequently we shall briefly discuss
the traditional approach to such problems.

Insurance as an institution gives its customers the ability to share
the risk they may face in the future by buying a suitable contract.
The law of averages or Strong Law of Large Numbers is used to reduce
risk by sharing a part of it between a large group of customers.
Given that $N$ individuals are willing to buy $N$ contracts of the
same type that pay a fixed amount $X_0=1$ if the defined risk, death
occurs during a time interval $[0,T]$ and by ignoring  fees and
taxes, the premium $p$ should be a function of $N,T,X_0$ and $q$---
the probability that the risk will happen during that interval. The
SLLN says then that if we have independence between different
individuals, then
$$
p=\E\left(e^{-\delta\,T}\,\frac{Y_N}{N}\right),
$$
where $\delta$ is the discount factor and $Y_N$ is the number of
customers who die, then $p=e^{-\delta\,T}\,q$ with
$$
e^{-\delta\,T}\,\frac{Y_N}{N}\rightarrow p\;\;a.s.,
$$
when $N$ goes to infinity, so that $p$ is a fair net premium to
charge for the insurance.

In case the size of the loss is uncertain then the premium is given
by $(1+\theta)\,p$, where $p$ is the premium for the average losses
and $\theta$ is a (safety) loading factor to cover possible
fluctuations.

In practice also, customers are of different ages so that $p$ varies
and it is assumed that the type of contract influences mortality
risk so that different values of $p$ are used for different types of
contracts.

In the presence of a financial risk (e.g equity-linked insurance
contracts), the direct application of the SLLN principle may not
give a suitable result as it does not take into account the
possibility of investing in the financial market and the restriction
of such pricing to purely financial claims does not necessarily
respect the no-arbitrage property.

As we can see, this procedure implies the use of a coherent risk
measure for valuing discounted monetary risks.

Many papers have been devoted to this kind of problem and many
techniques have been proposed to price such contracts. We recall the
risk-minimizing technique which considers the biometric risk as a
non-tradable risk in an incomplete financial market, see T. M\o ller
\cite{Moller} for more details.

In this paper, we propose to build a pricing that respects both SLLN
and no-arbitrage principles. In order to do this, we recall in
section 3 some results on one-period coherent risk measure and the
well-known theorem giving its representation in terms of test
probabilities. In section 4, we work in a multi-period case, we
define a chain of coherent risk measures that define prices along
the time axis and introduce some properties, namely lower, weak and
strong time-consistency. While lower time-consistency is a natural
property in this context, the weak one suggests that the pricing is
derived from a single set of test probabilities and the strong one
allows us to hedge a claim by a trade at each period. In section 5,
we consider the financial market as an embedded entity in the global
market and decompose a given pricing into its financial and
intermediary, or prerisk, parts. We show that the pricing can be
constructed from its two parts under the time-consistency property.
Finally, in section 7, we fix a no-arbitrage pricing mechanism $\Pi$
on the financial market and derive the family of time-consistent
pricing mechanisms that coincide with $\Pi$ on the purely financial
claims.

\section{One-period coherent risk measures.}
Let $(\Om,\calF,\P)$ be a probability space with $\calF_0\subset\calF$ a
sub-$\sigma$-algebra. In this section we recall the main result on
the characterization of a one-period coherent risk measure defined
on the vector space $\Linf(\calF)$ with values in $\Linf(\calF_0)$. The
$\sigma$-algebra $\calF_0$ is not necessarily trivial.

\begin{defi}
\label{d1}(See Delbaen \cite{delbaen}) We say that the mapping
$\rho:{\Linf}(\calF)\rightarrow {\Linf}(\calF_0)$ is a coherent risk
measure if it satisfies the following axioms :
\begin{enumerate}
\item Monotonicity: For every $X,Y\in {\Linf}(\calF)$,
$$
X\leq Y\,\mbox{a.s}\;\Rightarrow \;\rho(X)\leq \rho(Y)\;\mbox{a.s}.
$$
\item Subadditivity: For every $X,Y\in {\Linf}(\calF)$,
$$
\rho(X+Y)\leq\rho(X)+\rho(Y)\;\mbox{a.s}.
$$
\item Translation invariance: For every $X\in {\Linf}(\calF)$ and $y\in  {\Linf}(\calF_0)$,
$$\rho(X+y)=\rho(X)+y\;\mbox{a.s}.$$
\item $\calF_0$-Positive homogeneity: For every $X\in {\Linf}(\calF)$ and $a\in  {\Lplus}(\calF_0)$, we have
$$\rho(a\,X)=a\,\rho(X)\;\mbox{a.s}.$$
\end{enumerate}
\end{defi}

\begin{defi}
\label{d2} The coherent risk measure
$\rho:{\Linf}(\calF)\rightarrow {\Linf}(\calF_0)$ is said to
satisfy the Fatou property if a.s $\rho(X)\leq\liminf\rho(X_n)$,
for any sequence $(X_n)_{n\geq 1}$ uniformly bounded by $1$ and converging
to $X$ in probability.
\end{defi}

\begin{defi}
\label{d3}The coherent risk measure $\rho:{\Linf}(\calF)\rightarrow
{\Linf}(\calF_0)$ is called relevant if for each set $A\in\calF$ with
$\P[A|\, \calF_0]> 0$ a.s, we have that $\rho(1_A) > 0$ a.s.
\end{defi}

\begin{prop}
\label{p2}(See Delbaen \cite{delbaen}) Let the mapping
$\rho:{\Linf}(\calF)\rightarrow {\Linf}(\calF_0)$ be a relevant
coherent risk measure satisfying the Fatou property. Then
\begin{enumerate}
\item The acceptance set $\A_\rho:=\{X\in {\Linf}(\calF)\;;\;\rho(X)\leq 0\;\mbox{a.s}\}$
is a weak$^*$-closed convex cone, arbitrage-free, stable under
multiplication by bounded positive $\calF_0$-measurable random
variables and contains $\Linf_-(\calF)$.
\item There exists a
convex set of probability measures $\Q$, all of them being
absolutely continuous with respect to $\P$, with their densities forming an
$L^1(\P)$-closed set, and such that for $X\in {\Linf}(\calF)$:
\begin{eqnarray}
\label{ccmana}
\rho(X)=\mbox{ess-sup}\left\{\E_{\bbQ}(X|\,\calF_0)\;;\;\bbQ\in\Q
\right\}.
\end{eqnarray}
\end{enumerate}
\end{prop}
\begin{proof}We sketch the proof. Since $\rho$ is a coherent risk measure, $\A_\rho$ is a
convex cone, closed under multiplication by bounded positive
$\calF_0$-measurable random variables and contains $L^\infty_-$. Its
weak$^*$-closeness follows from the Fatou property and it is
arbitrage-free since $\rho$ is relevant. Now for the second
assertion, we remark that, by applying the Hahn-Banach separation
theorem with exhaustion argument (as in Schachermayer \cite{scha}),
we may deduce that there exists some $g\in \A_\rho^*$, where
$\A_\rho^*$ is the dual cone of $\A_\rho$ in $L^1$, such that $g>0$
a.s, then we define
$$
\Q ^e=\left\{\bbQ\ll\P\;;\;\frac{d\bbQ}{d\P}\in
\A_\rho^*,\frac{d\bbQ}{d\P}>0\;\mbox{a.s}\right\},
$$
and $\Q =\overline{\Q ^e}$. Now let $X\in \Linf(\calF)$ and
$f^+\in\Linf_+(\calF_0)$, then by the translation invariance property,
we get that $f^+\,(X-\rho(X))\in\A_\rho$ and for every
$\varepsilon>0$, $f^+\,(X-\rho(X))+\varepsilon\notin\A_\rho$. Consequently, we deduce (\ref{mana}).
\end{proof}

\begin{defi}
Given a coherent risk measure $\rho$, we define $\Q^\rho$ as
follows:
\begin{eqnarray}
\label{mana1} \Q^\rho =\left\{\bbQ\ll\P\;;\;\frac{d\bbQ}{d\P}\in
\A_\rho^*\right\}.
\end{eqnarray}
Conversely, given $\Q$ a collection (not necessarily closed, or
convex) of probability measures absolutely continuous with respect
to $\P$, we define
\begin{eqnarray}
\label{mana}
\rho^\Q(X)=\mbox{ess-sup}\left\{\E_{\bbQ}(X|\,\calF_0)\;;\;\bbQ\in\Q
\right\}.
\end{eqnarray}
The set $\Q^\rho$ is the largest subset $\Q$ for which
$\rho=\rho^\Q$.
\end{defi}

\section{Risk measure versus Market.}
Returning to our problem: we suppose that the intermediary is
equipped with the probability space $(\Om,\G,\P)$, with a filtration
$\bbG=(\G_t)_{t=0}^T$, with $\G=\G_T$, modelling the flow of
information on the discrete time axis $\bbT^+=\bbT\cup\{T\}$ with
$\bbT=\{0,...,T-1\}$.

We further suppose that the intermediary's pricing mechanism is
$\undrho=(\rho_0,...,\rho_{T-1})$ where each $\rho_t$ denotes the
price at time $t$ of future (discounted) payoffs. Note that by
choosing to price the discounted payoffs rather than the payoffs
themselves, it's not necessary to introduce the discount rate in the
property of translation invariance. Define the acceptance set of
positions
$$
\A^t=\{X\in
{\Linf}(\G)\;;\;\rho_t(X)\leq 0\,\, \P\hbox{ a.s. }\},
$$
the set of liabilities which the intermediary is willing to accept
for no nett charge or no nett reserve at time $t$.

\begin{defi}We say that the vector $\undrho=(\rho_0,...,\rho_{T-1})$ is a chain
of coherent risk measures if for each $t\in \bbT$, the mapping
$\rho_t:{\Linf}(\G_T)\to {\Linf}(\G_t)$ fulfills all the properties
of a relevant coherent risk measure with the Fatou property (taking
$\calF=\G_T, \calF_0=\G_t$ in Definition \ref{d1}).
\end{defi}

It follows from Proposition \ref{p2} that for all $t\in\bbT$, there
exists an $L^1$-closed convex set of probabilities $\Q
^t=\Q^{\rho_t}$, absolutely continuous w.r.t $\P$ such that for
every $X\in \Linf(\G)$,
$$
\rho_t(X)=\mbox{ess-sup}_{\bbQ\in\Q ^t}\E_\bbQ(X|\,\G_t).
$$
To determine the natural relationship between the subsets of
probability measures
$$(\Q ^0,...,\Q ^{T-1})=(\Q ^{\rho_0},...,\Q^{\rho_{T-1}}),
$$
let us consider a contract $C^{t,T}_X$ issued at time $t$ and paying
the $t$-discounted amount $X\in {\Linf}(\G)$ (i.e discounted to time
$t$) to the holder at time $T$. Its price at time $t$ is given by
$\rho_t(X)$. The buyer may choose, instead to buy another contract
$C^{t,t+s}_{\rho_{t+s}(X)}$ paying $\rho_{t+s}(X)$ at time $t+s$.
Its price is given by $\rho_t\circ\rho_{t+s}(X)$. This contract can
be seen also as a contract which gives the buyer, the right to
choose at time $t+s$ between cash $\rho_{t+s}(X)$ or a new contract
$C^{t+s,T}_X$. We conclude then that for every $t,t+s\in\bbT^+$ and
$X$ we should have
\begin{eqnarray}
\label{ltc}\rho_t(X)\leq\rho_t\circ\rho_{t+s}(X).
\end{eqnarray}
We say that $\undrho$ is lower time-consistent if $\undrho$
satisfies (\ref{ltc}) which is equivalent to saying that the
acceptance sets satisfy $\A^{t+s}\subset \A^t$ or by a duality
argument that $\Q ^t\subset\Q ^{t+s}$. In the case where the
inequality in (\ref{ltc}) becomes equality we say that $\undrho$ is
time-consistent w.r.t the filtration $(\G_t)_{t=0}^T$ or simply
$\bbG$-time-consistent.

\begin{defi}\label{def3}Let $t\geq s$ with $t,s\in \bbT$, $\calH$ and $\calH'$ be two
subsets of  probability measures on $(\Omega,\G)$. We say that
$\calH\subset_{s,t}\calH'$ if for every $\bbQ\in \calH$, there
exists some $\bbQ'\in \calH'$ such that for every $X\in
{\Linf}(\G_t)$, we have
$$
\E_\bbQ(X|\,\G_{s})=\E_{\bbQ'}(X|\,\G_{s}).
$$
We define the relation $\equiv_{s,t}$ in an analogous fashion and
$$
[\calH]_{s,t}=\left\{\bbQ\;\mbox{a  probability
measure}:\;\{\bbQ\}\subset_{s,t}\calH\right\}.
$$
For a $\P$-absolutely continuous probability measure $\R$,  we
denote by $\Lam^\R$ or $\Lam(\R)$ its density (with respect to $\P$)
and define $\Lam^\R_t=\E(\Lam^\R|\G_t)$ for every $t\in\bbT^+$, so
that $\Lam^\R_t$ is the density of the restriction of $\R$ to
$\G_t$.

\end{defi}

\begin{remark}The set $ [\calH]_{s,t}$ defined in the previous definition, is not necessarily closed in
$L^1$ even when $\calH$ is.
\end{remark}

\begin{defi}
Given a set of probability measures $\Q $,
\begin{enumerate}
\item
We define the associated chain of coherent risk measures
$\undrho^\Q=(\rho^\Q_0,\ldots,\rho^\Q_{T-1})$ as follows: for all
$t\in\bbT$ we define for $X\in\Linf(\G)$,
$$
\rho^\Q_t(X)=\mbox{ess-sup}\left\{\E_\bbQ(X|\G_t);\;\bbQ\in\Q\right\}.
$$
\item Let $\A $ be the dual cone of $\Q$; we define for $t\in\bbT$:
$$
\A_t=\{X;\;\al X\in \A\;\mbox{for all}\;\al\in \Linf_+(\G_t)\}.
$$
Remark that $\A=\A_0$ since $\G_0$ is trivial.
\end{enumerate}
\end{defi}

\begin{lem}\label{rt}
Given a set of probability measures $\Q $ with the dual cone $\A$,
then for all $t,t+s\in\bbT^+$, the dual cone of $[\Q]_{t,t+s}$ is
given by $\A_t\cap\Linf(\G_{t+s})+\Linf_-(\G_T)$.
\end{lem}
\begin{proof}Let $X\in\left([\Q]_{t,t+s}\right)^*$ and define $Y=X-\rho'_{t+s}(X)$ where
$\rho'$ is the associated coherent risk measure to the set
$[\Q]_{t,t+s}$, then $X=Y+\rho'_{t+s}(X)$. We want to show that
$Y\in\Linf_-(\G_T)$ and $\rho'_{t+s}(X)\in \A_t\cap\Linf(\G_{t+s})$.
Choose $g\in\Linf_+(\G_T)$ and define the probability $\bbQ$ having
the density
$$
f=\dfrac{g}{\E(g|\G_{t+s})}\Lam_{t+s},
$$
where $\Lam$ is the density of a probability measure $\R\in\Q$.
Remark that $\bbQ\equiv_{t,t+s}\R$, then $\bbQ\in[\Q]_{t,t+s}$ and
$$
\E g Y= \E_\bbQ \dfrac{\E(g|\G_{t+s})}{\Lam_{t+s}}
Y=\E_\bbQ\left(\dfrac{ \E(g|\G_{t+s})}{\Lam_{t+s}}\E_\bbQ
(Y|\G_{t+s})\right),
$$
with $\E_\bbQ (Y|\G_{t+s})=\E_\bbQ(X|\G_{t+s})-\rho'_{t+s}(X) \leq
0$. Hence $ \E g Y\leq 0$ for all $g\in\Linf_+(\G_T)$, which leads
to $Y\in\Linf_-(\G_T)$. Now Choose $\bbQ\in\Q$ and
$\al\in\Linf_+(\G_t)$, we have
$$
\E_\bbQ(\al \rho'_{t+s}(X))=\E (\Lam^\bbQ_{t+s}\al
\rho'_{t+s}(X))=a\E (\dfrac{\Lam^\bbQ_{t+s}\al}{a}
\rho'_{t+s}(X))=a\E (f\,\rho'_{t+s}(X)),
$$
with $a=\E(\Lam^\bbQ_t\al)$ and
$$
f=\dfrac{\Lam^\bbQ_{t+s}\al}{a}.
$$
Remark that there exists a sequence $\bbQ^n\in[\Q]_{t,t+s}$ such
that the increasing sequence $\E_{\bbQ^n}(X|\G_{t+s})$ converges a.s
to $\rho'_{t+s}(X)$. We denote by $\Lam^n$, the density of $\bbQ^n$.
We obtain
$$
\E_\bbQ(\al \rho'_{t+s}(X))=
a\lim_{n\rightarrow\infty}\E(f\dfrac{\Lam^n}{\Lam^n_{t+s}}X).
$$
Define
$$
f^n=f\dfrac{\Lam^n}{\Lam^n_{t+s}},
$$
and remark that $f^n_{t+s}=f^n_{t}=1$, then the associated
probability $\bbQ^n_1\in[\Q]_{t,t+s}$. In consequence
$$
\E_\bbQ(\al \rho'_{t+s}(X))\leq
a\lim_{n\rightarrow\infty}\E_{\bbQ^n_1}(X)\leq 0.
$$
Conversely let $\bbQ\in [\Q]_{t,t+s}$, there exists then some
$\bbQ'\in\Q$ such that $\bbQ\equiv_{t,t+s}\bbQ'$. We obtain for all
$X\in\A_t\cap\Linf(\G_{t+s})$,
$$
\E_\bbQ X=\E \Lam_{t+s}X=\E_{\bbQ'} \dfrac{\Lam_{t}}{\Lam'_{t}}X,
$$
with $\Lam$ and $\Lam'$ are respectively the densities of $\bbQ$ and
$\bbQ'$. Define $Z=\dfrac{\Lam_{t}}{\Lam'_{t}}$ and for all $n$,
define $Z^n=Z\1_{Z\leq n}$. Consequently
$$
\E_\bbQ X=\E ZX = \lim_{n\rightarrow \infty}\E_{\bbQ'} Z^n X\leq 0,
$$
since $Z^n X\in\A$ for all $n$.
\end{proof}

\begin{lem}\label{l22}Let $\Q$ be a set of $\P$-absolutely continuous probability measures on
$(\Omega,\G)$ with $\A$ its dual cone. Then
\begin{enumerate}
\item for all $t,s\in \bbT$,
$$(\A_t)_s=(\A_s)_t=\A_{t\vee s}.$$
\item for all $t,s\in \bbT^+$,
$$\overline{[\overline{[\Q]}_{s,T}]}_{t,T}=\overline{[\Q]}_{s\vee
t,T},
$$
with the closure taken in $L^1$.
\end{enumerate}
\end{lem}
\begin{proof}Suppose $s\leq t$, by definition $(\A_t)_s\subset\A_{t}$ and $(\A_s)_t\subset\A_{t}$, Now let
$X\in\A_{t}$, then $\al_t\,X\in \A$ for all $\al_t\in\Linf_+(\G_t)$
which means that $\be_s\,\al_t\,X\in \A$ for all
$\al_t\in\Linf_+(\G_t)$ and $\be_s\in\Linf_+(\G_s)$, we deduce that
$\al_t\,X\in \A_s$ (resp. $\be_s\,X\in \A_t$) for all
$\al_t\in\Linf_+(\G_t)$ (resp. for all $\be_s\in\Linf_+(\G_s)$),
therefore $X\in (\A_s)_t$ (resp. $X\in(\A_t)_s$). For the second
assertion, we apply Lemma \ref{rt} and the assertion $(1)$ and
obtain
$$
\left([\overline{[\Q]}_{s,T}]_{t,T}\right)^*
=\left(\A_{s}\right)_{t}=\A_{s\vee t}= \left([\Q]_{s\vee
t,T}\right)^*.
$$
\end{proof}

\begin{remark}
\label{x0} Given a set of probability measures $\Q $, the associated
chain of coherent risk measures $\undrho^\Q$ is lower
time-consistent.
\end{remark}

\begin{lem}
\label{l1} Let $\undrho=(\rho_0,...,\rho_{T-1})$ be a chain of
coherent risk measures with the associated vector of test
probabilities $(\Q ^0,...,\Q ^{T-1})$. Then there exists a single
$\Q $ such that for every $t\in \bbT$ we have $\rho_t=\rho^{\Q }_t$
on $\Linf(\G_{t+1})$. Moreover, if for every $t\in\{0,\ldots,T-2\}$
we have $\Q^{t+1}=\overline{[\Q^t]}_{t+1,T}$, then there exists a
single $\Q $ (that we can take equal to $\Q^0$ ) such that for every
$t\in \bbT$ we have $\rho_t=\rho^{\Q }_t$ on $\Linf(\G_{T})$.
\end{lem}
\begin{proof}Let us define the subset $ \Q
=\bigcap_{t=0}^{T-1}[\Q^t]_{t,t+1}$. Remark that for all $t\in\bbT$,
we have $\Q\subset[\Q^t]_{t,t+1}$, then
$[\Q]_{t,t+1}\subset[\Q^t]_{t,t+1}$. Now let
$\bbQ\in[\Q^t]_{t,t+1}$, then there exists some $\bbQ^t\in\Q^t$ such
that $\Q\equiv_{t,t+1}\bbQ^t$. Let $f^t$ denote the density of
$\bbQ$ and define $\bbQ'$ as the probability measure associated to
the density
$$
f=\prod\limits_{u\in\bbT}\frac{f^u_{u+1}}{f^u_u},
$$
where each $f^u$ is the density of a probability measure
$\bbQ^u\in\Q^u$ for $u\neq t$. Then $\bbQ\equiv_{t,t+1}\bbQ'$ with
$\bbQ'\equiv_{t,t+1}\bbQ^t\in\Q^t$ and for all $s\neq t$, we have
$\bbQ'\equiv_{s,s+1}\bbQ^s\in \Q^s$. In consequence $\bbQ\in
[\Q]_{t,t+1}$ and hence for all $t\in \bbT$ we have $\Q
\equiv_{t,t+1}\Q^t$. We deduce then that $\rho^\Q_t=\rho_t$ on
$\Linf(\G_{t+1})$.

Now suppose that for every $t\in\{0,\ldots,T-2\}$ we have
$\Q^{t+1}=\overline{[\Q^t]}_{t+1,T}$. We define $ \Q =\Q^0$ and
prove by induction on $t=1,\ldots,T-1$ that
$\Q^t=\overline{[\Q^0]}_{t,T}$. By assumption
$\Q^1=\overline{[\Q^0]}_{1,T}$, we suppose that the induction
hypothesis is true until $t$, then
$$
\Q^{t+1}=\overline{[\Q^t]}_{t+1,T}=\overline{[\overline{[\Q^0]}_{t,T}]}_{t+1,T}
=\overline{[\Q^0]}_{t+1,T},
$$
where the last equality is due to Lemma \ref{l22}.
\end{proof}

\begin{defi}
We say that a chain $\undrho=(\rho_0,...,\rho_{T-1})$ is weakly
time-consistent if there exists a single $\Q $ such that
$\undrho=\undrho^{\Q }$.
\end{defi}

\begin{cor}
\label{coronia} Let $\undrho=(\rho_0,...,\rho_{T-1})$ be a chain
with the associated vector of test probabilities $(\Q ^0,...,\Q
^{T-1})$. Then the chain is weakly time-consistent iff for every
$t\in\bbT$ we have $\Q ^{t}=\overline{[\Q ^0]}_{t,T}$.
\end{cor}

\begin{cor}
Let $\undrho=(\rho_0,...,\rho_{T-1})$ denote a chain of coherent
risk measures with the associated vector of test probabilities $(\Q
^0,...,\Q ^{T-1})$ and the family of dual cones
$(\A^0,\ldots,\A^{T-1})$ with $\A=\A^0$. Then $\undrho$ is weakly
time-consistent iff for all $t\in\bbT$ we have $\A^t=\A_t$.
\end{cor}

\begin{defi}
We say that a chain $\undrho=(\rho_0,...,\rho_{T-1})$ is
time-consistent if for every $s,t\in \bbT$ with $s\leq t$ we have
$\rho_s=\rho_s\circ\rho_t$.
\end{defi}

We note that every time-consistent chain is weakly time-consistent.
For the converse to hold, the maximal associated set $\Q $ of
probability measures has to satisfy the multiplicative stability
property (see Delbaen \cite{delbaen2}).

\begin{defi}We say that a set of $\P$-absolutely continuous probability measures $\Q $,
is $\bbG$-m-stable (or just m-stable if
there is no confusion as to the filtration) if for every $\bbQ\in\Q
$, $\bbQ'\in\Q ^e$ and $t\in \bbT$, the probability measure ${\tilde
\bbQ}$ is contained in $\Q $, where
$$
\Lam^{\tilde \bbQ}=
\Lambda^\bbQ_t\;\dfrac{\Lambda^{\bbQ'}}{\Lambda^{\bbQ'}_t}.
$$
\end{defi}

\begin{remark}
The property of m-stability was defined by  Delbaen in
\cite{delbaen2} and the property was first introduced for EMMs by
Jacka in \cite{jacka}.
\end{remark}

\begin{lem}
\label{78}Let $\Q$ be a set of $\P$-absolutely continuous
probability measures, then $\bigcap_{t\in\bbT}[\Q]_{t,t+1}$ is the
smallest m-stable set of probability measures containing $\Q$ and
therefore $\Q$ is m-stable iff $\Q=\bigcap_{t\in\bbT}[\Q]_{t,t+1}$.
\end{lem}
\begin{proof}Let us define $\calH\defto\bigcap_{t\in\bbT}[\Q]_{t,t+1}$. We show first
that $\calH$ is m-stable. In order to do this, let $t\in \bbT$,
$\bbQ\in\calH$ and $\bbQ'\in \calH^e $ with respective densities
$\Lam$ and $\Lam'$. Define the probability measure $\R$ by
$$
\Lam(\R)=\Lam_t\;\dfrac{\Lam'}{\Lam'_t}.
$$
We want to show that $\R\in \calH $, so it remains to show that
$\R\in [\Q]_{s,s+1}$ for all $s\in\bbT$. Remark that
$$
\dfrac{\Lam_{s+1}(\R)}{\Lam_s(\R)}=\left\{
\begin{array}{l}
\dfrac{\Lam'_{s+1}}{\Lam'_s} \;\;\;\;\mbox{for}\; \;s\geq t\\
\\
\dfrac{\Lam_{s+1}}{\Lam_s} \;\;\;\;\mbox{for} \;\;s\leq t-1.
\end{array}\right.
$$
In consequence $\R\equiv_{s,s+1}\bbQ'$ for $s\geq t$ and
$\R\equiv_{s,s+1}\bbQ$ for $s\leq t-1$. Remark that since
$\bbQ,\bbQ'\in\calH\subset[\Q]_{s,s+1}$, then there exists
$\bbQ_s,\bbQ'_s\in\Q$ such that $\bbQ'\equiv_{s,s+1}\bbQ'_s$ and
$\bbQ\equiv_{s,s+1}\bbQ_s$, therefore $\R\equiv_{s,s+1}\bbQ'_s$ for
$s\geq t$ and $\R\equiv_{s,s+1}\bbQ_s$ for $s\leq t-1$ with
$\bbQ_s,\bbQ'_s\in\Q$.

Now let $\calH'$ be an m-stable set of probability measures
containing $\Q$ and let $\bbQ\in\calH$, then there exists
$\bbQ^0,\ldots,\bbQ^{T-1}\in\Q$ with their respective densities
$\Lam^0,\ldots,\Lam^{T-1}$ such that
$$
\Lam^\bbQ=\prod\limits_{u\in\bbT}\dfrac{\Lam^{u}_{u+1}}{\Lam^{u}_{u}}.
$$
We define, for each $t\in\bbT$,
$$
Z^t=\Lam^{t}_{t+1}\;\prod\limits_{u=t+1}^{T-1}\dfrac{\Lam^{u}_{u+1}}{\Lam^{u}_{u}}.
$$
Remark that $\Lam^\bbQ=Z^0$, then we prove by induction on
$t=T-1,\ldots,0$ that $Z^t\in\calH'$. We have $Z^{T-1}=\Lam^{T-1}\in
\Q\subset\calH'$, now suppose that $Z^{t+1}\in\calH'$ and remark
that
$$
Z^t=\dfrac{Z^{t+1}}{Z^{t+1}_{t+1}}\,\Lam^{t}_{t+1}.
$$
Since $\calH'$ is m-stable, we obtain $Z^t\in\calH'$. The
equivalence in Lemma \ref{78} becomes straightforward.
\end{proof}

The following theorem is due to Delbaen (\cite{delbaen2}).
\begin{theorem}Given a set of probability measures $\Q $ (not necessarily a
closed convex set), the associated chain $\undrho^\Q $ is
time-consistent iff $\overline{co(\Q )}$ is m-stable. By a small
abuse of language we say that $\Q $ is time-consistent when
$\undrho^\Q $ is.
\end{theorem}

Here we state some simple and interesting results on
time-consistency.
\begin{theorem}\label{domin}
\label{t-c2} Let $\undrho=(\rho_0,...,\rho_{T-1})$ be a lower
time-consistent chain and define for each $t\in \bbT$, the risk
measure $\eta_t=\rho_t\circ...\circ\rho_{T-1}$. Then ${\underline
\eta}\defto(\eta_0,...,\eta_{T-1})$ is the minimal time-consistent
chain which dominates $\undrho$.
\end{theorem}
\begin{proof}By definition $\eta_t=\rho_t\circ\eta_{t+1}$ and $\eta_t=\rho_t$ on
${\Linf}(\G_{t+1})$, so ${\underline \eta}$ is time-consistent. The
fact that $\eta$ dominates $\rho$ follows by backwards induction.

Now let ${\underline \xi}\defto(\xi_0,...,\xi_{T-1})$ be a
time-consistent chain of coherent risk measures which dominates
$\undrho$. Therefore $\xi_{T-1}\geq \rho_{T-1}=\eta_{T-1}$ and by
backwards induction on $t$ we have
$$
\xi_t=\xi_t\circ\xi_{t+1}\geq
\xi_t\circ\eta_{t+1}\geq\rho_t\circ\eta_{t+1}=\eta_t.
$$
\end{proof}

\begin{remark}$\eta$ corresponds to the smallest m-stable, closed convex set
of probability measures containing $\Q^\rho$, i.e
$\Q^\eta=\overline{\bigcap_{t\in\bbT}[\Q^\rho]_{t,t+1}}$.
\end{remark}

\begin{lem}
\label{t-c}Suppose that $\Q $ is time-consistent and $s\in \bbT$.
Let $X\in {\Linf}(\G_T)$ be such that there exists a probability
measure $\bbQ\in \Q^e $ satisfying $\rho_s(X)=\E_{\bbQ}(X|\G_s)$.
Then for every $t\geq s$ we have:
$$
\rho_{t}(X)=\E_{\bbQ}(X|\G_{t}).
$$
\end{lem}

\begin{proof} It suffices to remark that $\rho_{t}(X)\geq
\E_{\bbQ}(X|\G_t)$ and
$$
\E_{\bbQ}\rho_{t}(X) = \E_{\bbQ}\E_{\bbQ}(\rho_{t}(X)|\G_s)
\leq\E_{\bbQ}\rho_s\circ\rho_{t}(X)=\E_{\bbQ}\rho_{s}(X).
$$
It's given that $\rho_s(X)=\E_{\bbQ}(X|\G_s)$, then
$$
\E_{\bbQ}\rho_{t}(X)\leq
\E_{\bbQ}\E_{\bbQ}(X|\G_s)=\E_{\bbQ}\,X=\E_{\bbQ}\E_{\bbQ}(X|\G_t).
$$
Then
$$
\E_{\bbQ}\left(\rho_t(X)-\E_{\bbQ}(X|\G_t)\right)=0,
$$
which means that a.s
$$
\rho_t(X)=\E_{\bbQ}(X|\G_t).
$$
\end{proof}

\begin{theorem}
\label{t-cbab}$\Q$ is time-consistent iff the process
$(\rho^\Q_t(X))_{0\leq t\leq T}$ is a $\Q$-uniform-supermartingale
for every $X\in {\Linf}(\G_T)$.
\end{theorem}

\begin{proof}Suppose that $\Q $ is time-consistent. Then for $X\in
{\Linf}(\G_T)$ and $\bbQ\in \Q $ we have (suppressing the $\Q$-dependence of $\rho$):

$$
\E_\bbQ(\rho_{t+s}(X)|\G_t)\leq
\rho_{t}\circ\rho_{t+s}(X)=\rho_{t}(X).
$$
Now suppose that the process $(\rho_t(X))_{t=0}^T$ is a
$\Q$-uniform-supermartingale; which means that for every $\bbQ\in \Q $:
$$
\E_\bbQ(\rho_{t+s}(X)|\G_t)\leq \rho_{t}(X).
$$
It follows that $\rho_{t}\circ\rho_{t+s}(X)\leq \rho_{t}(X)$ and the
result follows from lower time-consistency.
\end{proof}

In the next result we show the relationship between the
time-consistency of the chain $\undrho$ and the decomposition of its
acceptance set $\A=\{X\in\Linf\;;\;\rho_0(X)\leq 0\}$. Define for
every $t\in \bbT$,
$$\K_t\defto \{X\in \Linf(\G_{t+1})\;;\;\rho_t(X)\leq
0\}=\A_t\cap\Linf(\G_{t+1}).
$$

\begin{theorem}
\label{bara}Suppose that the chain $\undrho$ is lower
time-consistent, then it is time-consistent iff
$\A=\K_0+...+\K_{T-1}$. In this case for all $t\in\bbT$, we have
$\A_t=\K_t+...+\K_{T-1}$.

\end{theorem}
\begin{proof}Suppose that $\undrho$ is time-consistent, then for every $X\in\A$ we get
$$
X=\sum_{s=1}^{T-1}\,(\rho_{s+1}(X)-\rho_s(X))+\rho_1(X),
$$
with $\rho_T=id$. Let $u_s=\rho_{s+1}(X)-\rho_s(X)\in \K_s$ for
$s\in\{1,...,T-1\}$ and $u_0=\rho_{1}(X)\in \K_0$. It follows that
$\A\subset \K_0+...+\K_{T-1}$. Since $\K_0+...+\K_{T-1}\subset\A$ we
have equality.

Now suppose that $\A=\K_0+...+\K_{T-1}$. Let $X\in\Linf$ and $t\in
\bbT$ be fixed. Since $\rho_t(X-\rho_t(X))=0$ and
$\rho_0\circ\rho_t\geq \rho_0$ it follows that $X-\rho_t(X)\in\A$,
and so there exist $y_0\in \K_0,...,y_{T-1}\in \K_{T-1}$ such that
$$
X-\rho_t(X)=y_0+...+y_{T-1}. $$ By applying $\rho_t$ to both sides of
this equality,
we obtain
$$
0=y_0+...+y_{t-1}+\rho_t(y_t+...+y_{T-1}),$$ and so by subadditivity
$$
0\leq y_0+...+y_{t-1}+\sum_{s=t}^{T-1}\,\rho_t(y_s)$$ and by lower
time-consistency and the assumption that $y_s\in \K_s$, we get
$$
0\leq y_0+...+y_{t-1}+\sum_{s=t}^{T-1}\,\rho_t\circ\rho_s(y_s)\leq
y_0+...+y_{t-1}.
$$
But $y_0+...+y_{t-1}\in\A$, so $\E_{\bbQ}(y_0+...+y_{t-1})\leq 0$ for
some $\bbQ\in\Q^e$ and therefore $y_0+...+y_{t-1}=0$.

Now it follows that $X-\rho_t(X)=y_t+...+y_{T-1}$. By successively
applying $\rho_{T-1},\ldots,\rho_t$ on both sides and using the
properties of $\rho_s$ and $\K_s$, we obtain
$$
\rho_t\circ\ldots\circ\rho_{T-1}(X)-\rho_t(X)=\eta_t(y_t+...+y_{T-1})\leq
0.
$$
Finally, since $\rho_t\circ\ldots\circ\rho_{T-1}(X)\geq \rho_t(X)$
it follows that $\rho_t\circ\ldots\circ\rho_{T-1}(X)=\rho_t(X)$ and
hence, by Lemma \ref{domin}, it follows that $\undrho$ is
time-consistent.
\end{proof}

\begin{remark}
It is in this situation (where $\undrho$ is time-consistent) that we
can replicate claims in $\A$ by a sequence of one-period trades.
This explains the \lq mark-to-market' requirement of section $1$.
\end{remark}

\section{The decomposition of the global market.}
\begin{Ex}
\label{ex1} We consider a contract that provides one share of $XYZ$
stock to the insured if he or she is still alive in one year's time,
and nothing otherwise.
\begin{enumerate}
\item What's the fair premium for this contract?
\item What's the \lq self financing ' strategy if it exists?
\end{enumerate}
To formulate this problem, let $S$ denote the discounted price of
the $XYZ$ share in one year's time and let
$$
Y=\left\{\begin{array}{ll} 1&\mbox{if the insured is alive then}\\
0&\mbox{otherwise}
\end{array}\right.
$$
We suppose that $S$ and $Y$ are defined respectively on two
probability spaces $(\Om^1,\G^1,\P^1)$ and $(\Om^2,\G^2,\P^2)$, with
$\G^1=\sigma(S)$ and that the pricing of purely financial (resp.
insurance) claims is given by $p_F$ (resp. $p_I$). The payoff of the
contract is $H=S\,Y$. To price such a claim in one-period case, we
remark first that $H=H(S)$ where $H(x)\defto x\,Y$ for a scalar $x$.
Remark also that for a fixed $x$, the claim $x\,Y$ is a purely
insurance claim and it's priced by $p_I(x\,Y)$ and that the claim
$H^F\defto p_I(x\,Y)|_{x=S}$ is a purely financial claim, priced by
$p_F(H^F)$. We propose then the premium of $H$, $p(H)=p_F(H^F)$. The
\lq self financing ' strategy will be the one to hedge the claim
$H^F$. We obtain then the decomposition of the claim $H$ as follows:
$H=p(H)+U^F+U^I$, where the claims $U^F\defto H^F-p_F(H^F)$ and
$U^I\defto H-H^F$ are admissible.
\end{Ex}

Under the assumption that both $P_F$ and $P_I$ are coherent risk
measures with $\A^F$ and $\A^I$ their respective acceptance sets,
the risk measure $p$ defined above is a coherent risk measure with
acceptance set $\A$, satisfying $\A=\A^F+\A^I_s$, where
$$
\A^I_s=\left\{X:\;f(S)X\in\A^I,\;\mbox{for
all}\;f\in\Linf_+(\R)\right\}.
$$

In this section, we suppose that we're given a probability space
$(\Om,\G,\P)$ and a coherent risk measure $\rho$ and $\F\subset\G$
the financial sub-$\sigma$-algebra. Our aim is
\begin{enumerate}\item to construct, first in the one-period case, two
coherent risk measures
$$
\rho_F:\Linf(\F)\rightarrow\R
\;\mbox{and}\;\rho_I:\Linf(\G)\rightarrow\Linf(\F),
$$
such that
$\rho=\rho_F$ on $\Linf(\F)$ and conditioning on $\F$,
$\rho=\rho_I$.
\item to establish necessary and sufficient conditions on $\rho$
such that $\rho=\rho_F\circ\rho_I$.
\end{enumerate}
Remark that if we denote
$\G_0=\{\emptyset,\Om\},\,\G_{0^+}=\F,\,\G_1=\G$ and suppose that
$\rho$ is time-consistent w.r.t the filtration
$(\G_0,\G_{0^+},\G_1)$, then the acceptance set $\A$ will be
decomposed as follows: $\A=\A^F+\A^I$ where $\A^F=\A\cap\Linf(\F)$
is the financial part of the whole market, whilst the second
component $\A^I=\A_{0^+}$, is the intermediary market, equivalent to
the whole market in the absence of the financial market
$(\calF=\G_0)$. Any claim then can be decomposed into its financial
and intermediary parts. The corresponding pricing mechanisms are
given respectively by $\Q^F=[\Q^\rho]_{0,0^+}$ and
$\Q^I=[\Q^\rho]_{0^+,1}$. Here we adopted the notation of the last
section.

To generalize this setting to a multi-period case, we introduce the
following notation. Let $(\Om,\G,\P)$ be a probability space
equipped with the filtration $(\G_t)_{t\in\bbT^+}$. Let
$(\calF_t)_{t\in\bbT^+}$ be the filtration modeling the information
in the financial market such that for every $t\in \bbT^+$ we have
$\calF_t\subset\G_t$ with $\calF_0$ and $\G_0$ trivial. We assume
that the intermediary makes prices according to a pricing mechanism
$\rho$, defined by a set of $\P$-absolutely continuous probabilities
$\Q $ on $\Om$. We suppose w.l.o.g that $\P\in \Q $ and that the set
$\Q $ is an $L^1(\P)$-closed convex set. Define $\Q ^e$ to be the
subset of $\P$-equivalent probability measures in $\Q $, the
intermediate $\sigma$-algebras $\G_{t^+}=\G_t\bigvee\calF_{t+1}$ and
the filtration
$$
\bbG^*=(\G_0,\G_{0^+},\G_1,...,\G_{T}).
$$

Define the subsets $\Q ^F$ and $\Q ^I$ as follows. For $t\in \bbT$,
we define $\Q ^{t,F}=[\Q]_{t,t^+}$ and $\Q ^F=\bigcap_{t=0}^{T-1}\Q
^{t,F}$. In the same way we define $\Q ^{t,I}=[\Q]_{t^+,t+1}$ and
$\Q ^I=\bigcap_{t=0}^{T-1}\Q ^{t,I}$. We denote respectively by
$\rho^F$ and $\rho^I$ the coherent risk measures associated to the
subsets $\Q ^F$ and $\Q ^I$. In the following lemmas we state some
interesting properties of these two subsets of probabilities.

\begin{defi}Let $t\in \bbT$. We define the binary relation $\sim_{t,F}$, defined on the
set of all $\P$-absolutely continuous probabilities,  as follows:
$$\bbQ\sim_{t,F}\bbQ' \hbox{ iff }\{\bbQ\}^{t,F}=\{\bbQ'\}^{t,F}.$$

We define $\sim_{t,I}$ in the same fashion.
\end{defi}
\begin{lem}$\sim_{t,F}$ is an equivalence
relation. Moreover $\Q ^{t,F}=\bigcup_{\bbQ\in\Q }\{\bbQ\}^{t,F}$.
The analogous results hold for $\sim_{t,I}$.
\end{lem}
\begin{proof}The binary relation $\sim_{t,F}$ is obviously an equivalence relation. Take
a probability measures $\bbQ\in \Q^{t,F}$, then there exists a
probability measure $\bbQ'\in\Q$ such that $\bbQ\sim_{t,F}\bbQ'$,
which means that $\bbQ\in\{\bbQ'\}^{t,F}$ and hence
$\Q^{t,F}\subset\bigcap_{\bbQ\in\Q}\{\bbQ\}^{t,F}$. The reverse
inclusion is obvious.
\end{proof}

\begin{lem}
\label{bml4} For every $\bbQ\in\Q ^F$ and $t\in \bbT$, there exists
some $\bbQ^t\in \Q $ such that $\E_\bbQ(X|\G_t)=\E_{\bbQ^t}(X|\G_t)$
for every $X\in {\Linf}(\G_{t^+})$. Analogously, for every $\bbQ\in
\Q ^I$ and $t\in \bbT$, there exists some $\bbQ^{t^+}\in \Q $ such
that $\E_\bbQ(X|\G_{t^+})=\E_{\bbQ^{t^+}}(X|\G_{t^+})$ for every
$X\in {\Linf}(\G_{t+1})$.
\end{lem}
\begin{proof}Immediate consequence of Definition \ref{def3}.
\end{proof}

\begin{theorem}
\label{bml4kk}Let $\Q $ be a set of $\P$-absolutely continuous
probability measures on $\Om$. Then
\begin{enumerate}
\item $\Q $ is $\bbG^*$-m-stable if and only if
$\Q =\Q ^F\cap\Q ^I$.
\item The subsets $\Q ^F$ and $\Q ^I$ are $\bbG^*$-m-stable. Moreover
$$\overline{(\Q ^F)^I}=\overline{(\Q ^I)^F}=\calP,
$$
where $\calP$ is the set of all $\P$-absolutely continuous
probability measures and the closure is taken in $L^1(\Om)$.
\end{enumerate}
\end{theorem}
\begin{proof}The first assertion is an immediate consequence of
Lemma \ref{78}. The second assertion is an immediate consequence of
assertion $(1)$ since $(\Q ^F)^F=\Q ^F$ and $\Q ^F\subset(\Q ^F)^I$,
so $\Q ^F=(\Q ^F)^F\cap(\Q ^F)^I$. We make the same argument for the
$I$-part. We remark also that $\calP^e\subset (\Q ^F)^I,(\Q ^I)^F$
where $\calP^e$ is the set of all $\P$-equivalent probability
measures. Indeed let $\bbQ\in\calP^e$ with $f=\Lam^\bbQ$ and $t\in
\bbT$ fixed. We define the probability $\bbQ^t$ by its density
$$
\Lam^{t}=\dfrac{f_{t+1}}{f_{t^+}}.
$$
Then $\bbQ^t\in\Q ^F$ since
$\bbQ^t\sim_{s,F}\P$ for every $s\in \bbT$. Moreover $\bbQ
\sim_{t,I} \bbQ^t$, therefore $\bbQ\in (\Q ^F)^I$. We do the same
for the inclusion $\calP^e\subset(\Q ^I)^F$.
\end{proof}

Let $\undrho=\undrho^\Q $ be the chain associated to the set of
probabilities $\Q $ and $\rho=\rho_0$. Let us define the acceptance
cone $\A=\A_{\rho}$ associated to the coherent risk measure $\rho$
by
$$
\A=\{X\in {\Linf}(\G)\;;\;\rho(X)\leq 0\}.
$$
$\A$ is then a weak$^*$-closed convex cone in $\Linf$. Our objective
is to decompose this trading cone in the global market into the sum
of two trading cones, one in the financial market and the other in
the intermediary's market.

We define the following convex cones
$$
\K_t^F=\{X\in {\Linf}(\G_{t^+})\;;\;\rho_t(X)\leq
0\}=\A_t\cap\Linf(\G_{t^+}),
$$
$$ \K_t^I=\{X\in
{\Linf}(\G_{t+1})\;;\;\rho_{t^+}(X)\leq
0\}=\A_{t^+}\cap\Linf(\G_{t+1}),
$$
for $t\in \bbT$, $\K^F_T=\Linf_-(\G)$,
$$
\A^F=\K_0^F+...+\K_{T}^F,
$$
and
$$
\A^I=\K_0^I+...+\K_{T-1}^I.
$$

Then

\begin{lem}
\label{bml2} $\A=\A^F+\A^I$ iff $\Q $ is $\bbG^*$-time-consistent.
\end{lem}

\begin{proof}Immediate consequence of Lemma
\ref{bara}.
\end{proof}
\begin{remark}Remark that this corresponds to mark-to-market approach valuation.
\end{remark}

The question now is to characterize the pricing mechanism in both
trading cones. In the next lemma we prove that the cones $\A^F$ and
$\A^I$ are respectively the acceptance sets of the risk measures
$\rho^F$ and $\rho^I$.

\begin{lem}
\label{bml3} $\A^F=\A_{\rho^F}$ and $\A^I=\A_{\rho^I}$.
\end{lem}

\begin{proof}Since $\Q^F$ is $\bbG^*$-time-consistent, then $\A_{\rho^F}=W+M$
where $W=W_0+...+W_{T}$ and $M=M_0+...+M_{T-1}$ with
$$
W_t=\{X\in {\Linf}(\G_{t^+})\;;\;\rho^F_t(X)\leq 0\},
$$
and
$$
M_t=\{X\in {\Linf}(\G_{t+1})\;;\;\rho^F_{t^+}(X)\leq 0\},
$$
for $t\in \bbT$ and $W_T=\Linf_-(\G)$. By definition $\Q \subset \Q
^F$, and we deduce that each $W_t\subset \K^F_t$. Now let $X\in
\K^F_t$ and $\bbQ\in \Q ^F$. By applying Lemma \ref{bml4}, there
exists some $\bbQ^t\in \Q $ such that
$$
\E_\bbQ(X|\G_t)=\E_{\bbQ^t}(X|\G_t)\leq \rho_t(X)\leq 0.
$$
In consequence $\rho^F_t(X)\leq 0$ and $X\in W_t$. We have that
$\K^F_t=W_t$ and in consequence $\A^F=W$. It suffices therefore to
prove that $M\subset \Linf_-$ which follows if we can prove that
each $M_t\subset \Linf_-$ for $t\in \bbT$.

Let $X\in M_t$, which means that $X\in {\Linf}(\G_{t+1})$ and
$\rho^F_{t^+}(X)\leq 0$. Then for every $\bbQ\in \Q ^F$ we have
$\E_\bbQ(X|\G_{t^+})\leq 0$. Now let $g\in \Linf_+(\G_{t+1})$ with
$g>0$ a.s and $\E g=1$. Define the probability measure $\bbQ_g$ by
$d\bbQ_g=fd\P$ where
$$
f=\dfrac{g}{\E(g|\G_{t^+})}. $$ Then $\bbQ_g\ll \P$ and for every
$s\in \bbT$ and $B\in \G_{s^+}$ we have
$$
\bbQ_g(B|\G_s)=\E\left(B,\dfrac{f_{s^+}}{f_s}|\G_s\right)=\P(B|\G_s),$$
since
$$
\dfrac{f_{s^+}}{f_s}:=\dfrac{\E(f|\G_{s^+})}{\E(f|\G_s)}=1.
$$
Consequently $\bbQ_g\in \Q ^F$ and so
$$
\E(g.X)=\E(\E(g|\G_{t^+})\,f\,X)=\E\left(\E(g|\G_{t^+})\,\E_{\bbQ_g}(X|\G_{t^+})\right)\leq
0,
$$
for every $g\in L^{1}_{+}(\G_{t+1})$. Therefore $X\leq 0$. In
the same way we prove that $\A^I=\A_{\rho^I}$.
\end{proof}

\begin{cor}
\label{bmc1} The convex cones $\A^F$ and $\A^I$ are weak$^*$-closed
in $\Linf$.
\end{cor}

\begin{cor}
\label{cad} Let $\Q_1$ and $\Q_2$ be two convex subsets in $\calP$
with $\rho^1$ and $\rho^2$ their respective coherent risk measures.
Then the following assertions are equivalent.
\begin{enumerate}
\item For all $t\in \bbT$, $\rho^2_t\leq\rho^1_t$ on
$\Linf(\G_{t^+})$.
\item $\A^{\Q^F_1}\subset\A^{\Q^F_2}$.
\item $\overline{\Q^F_2}\subset\overline{\Q^F_1}$.
\end{enumerate}
\end{cor}

\begin{proof}The assertions $(2)$ and $(3)$ are equivalent by duality argument. Now let suppose that
$(2)$ is satisfied, then for all $t\in\bbT$, we have
$$
K^1_t\defto
\A^{\Q_1}_t\cap\Linf(\G_{t^+})=\A^{\Q^F_1}_t\cap\Linf(\G_{t^+})\defto
K^{1,F}_t\subset \A^{\Q^F_2}_t\cap\Linf(\G_{t^+})=K^{2,F}_t=K^2_t.
$$
Take $X\in\Linf(\G_{t^+})$, then $X-\rho^1_t(X)\in K^1_t\subset
K^2_t$. Therefore $\rho^2_t(X-\rho^1_t(X))\leq 0$ which means that
$\rho^2_t(X)\leq \rho^1_t(X)$. Conversely for all $t\in\bbT$,
$$
K^{1,F}_t= K^{1}_t\subset K^{2}_t=K^{2,F}_t.
$$
The assertion $(2)$ is obtained.
\end{proof}

\section{Example.}
Consider the example, where sample spaces are $I=\{i,i'\}$,
$F=\{f,f'\}$, $T=1$ and the probabilities $\bbI$ and $\bbF$ are
given by $\bbI(i)=\bbF(f)=1/2$. The financial market can be seen
then as associated to one risky asset taking only two values at time
$1$ and a constant interest rate. This market is complete and we
suppose that $\bbF$ is the equivalent martingale measure. The sample
space is given by $\Om=I\times F=\{(i,f),(i,f'),(i',f),(i',f')\}$,
the probability measure $\P=\bbI\otimes\bbF$ and ${\Linf}(\Om)$ is
identified with the space of $2\times 2$-matrices. We define the
pricing set $\Q $ by:
$$\Q =\left\{\bbQ\ll\P\;;\;\Lambda^\bbQ\leq 1+\varepsilon\;\mbox{and}\;\Lambda^\bbQ_{0^+}=1\right\}.$$
The subset $\Q $ can also be written as follows
$$\Q =\left\{(q_{ij})_{1\leq i,j\leq 2}\;;\;\sum_{i,j}q_{ij}=1,\;0\leq q_{ij}\leq
1/4(1+\varepsilon) \;\mbox{and for
each}\;j:\;q_{1j}+q_{2j}=1/2\right\}.$$ Note that $\Q$ is chosen to
have margin $\bbF$ on $F$, and to correspond to a TailVaR type
construction on $I$.

To compute the corresponding quantities $\rho_{0^+}(X)$ for
$X\in\Linf(\Om)$, we remark that the extreme points of the set $\Q $
are given by
$$
\bbQ^{a,b}= \dfrac{1}{4}\left(\begin{array}{rl}
1+a\vare & 1+b\vare \\
&\\
1-a\vare & 1-b\vare
\end{array}\right)
$$
with $(a,b)\in\{(1,1),(1,-1),(-1,1),(-1,-1)\}$. Therefore we may
check easily that for $\om\in F$
$$\rho_{0^+}(\1_{(i,\om)})=\dfrac{1}{2}(1+\vare)\1_{(\om)},$$
and
$$\rho_{0^+}(-\1_{(i,\om)})=-\dfrac{1}{2}(1-\vare)\1_{(\om)}.$$
That means that for a real $x$ we have:
$$
\rho_{0^+}(x\1_{(i,f)})= \dfrac{1}{2}(x+\vare|x|)\1_{(f)}.
$$
Consequently, for every $X\in {\Linf}(\G_{1})$ we have:
$$
\rho_{0^+}(X)=\al_X(f)\1_{f}+\al_X(f')\1_{f'},
$$
with
$$
\al_X(g)=\dfrac{1}{2}\left(X(i,g)+\vare|X(i,g)|+X(i',g)+\vare|X(i',g)|\right).
$$
For every $X\in {\Linf}(\G_{0^+})$ we have:
$$
\rho_{0}(X)= \E(X).
$$
We conclude then that
$$
\Q ^I=\left\{(q_{ij})_{1\leq i,j\leq
2}\;;\;\sum_{i,j}q_{ij}=1,q_{ij}\geq 0\;\mbox{and for
each}\;j:\;\dfrac{1}{\delta_\vare}\leq \dfrac{q_{1j}}{q_{2j}}\leq
\delta_\vare\right\}.
$$
with $\delta_\vare=\dfrac{1+\vare}{1-\vare}$. Moreover
$$
\Q ^F=\left\{\left(\begin{array}{cc}
\al&\beta\\
1/2-\al&1/2-\beta
\end{array}\right)\;;\;0\leq \al,\beta\leq 1/2\right\}.
$$
The set $\Q $ is $\bbG^*$-time-consistent since
$\Q =\Q ^F\cap\Q ^I$.

\section{Pricing.}
In this section, we suppose we are in the same situation as in
Example \ref{ex1}, where the financial market is equipped with a
no-arbitrage pricing $\Pi$ (namely a closed convex set of
probability measures); defined on a probability space
$(\Om_F,\F,\P_F)$ with the filtration $(\F_t)_{t\in\bbT^+}$.
Moreover we consider a probability space $(\Om_I,\calI,\P_I)$ with
the filtration $(\calI_t)_{t\in\bbT^+}$, to model the biometric
risk.

Our aim is to build the class of pricing mechanism $\rho$ (or
$\Q^\rho$) that prices the purely financial claims as $\Pi$ does.

Define the product probability space $(\Om,\G,\P)$ as follows:
$\Om=\Om_F\times\Om_I,\G=\F\otimes\calI$ and $\P=\P_F\otimes\P_I$,
equipped with the filtration $(\G_t)_{t\in\bbT^+}$ given by
$\G_t=\F_t\otimes\calI_t$ and $\G_{t^+}=\F_{t+1}\otimes\calI_t$. Let
$\hatpi$ denote the extension of $\Pi$ to the product space, i.e $$
\hatpi=\{\bbQ\otimes\P_I:\;\bbQ\in\Pi\}. $$

We state first the following result and identify the probability
measures with their densities.
\begin{lem}
\label{ca}Let $\Q_1,\Q_2\subset\calP$, then the set $\Q$ defined by
$\Q =\Q_1^F\cap\Q_2^I$, satisfies the following:
$\Q^F=\Q_1^F$,\,$\Q^I=\Q_2^I$ and then $\Q$ is $\bbG^*$-m-stable.
\end{lem}

\begin{proof}
Remark that $\Q^F\subset(\Q^F_1)^F=\Q^F_1$ since $\Q\subset\Q^F_1$.
Now let $\bbQ\in\Q^F_1$, then for all $t\in\bbT$, there exists a
probability measure $\R^t\in\Q_1$ such that $\bbQ\sim_{t,F}\R^t$.
Define the probability measure $\bbQ^t$ by the density
$$
f^t=\prod\limits_{u\in\bbT}\left(\dfrac{\Lam^u_{u^+}}{\Lam^u_u}\,\dfrac{\Lam^{u^+}_{u+1}}{\Lam^{u^+}_{u^+}}
\right),
$$
where $\Lam^u$ and $\Lam^{u^+}$ are respectively the densities of
probability measures $\R^u\in\Q_1$ and $\R^{u^+}\in\Q_2$ for all
$u\in\bbT$. Remark that $\bbQ\sim_{t,F}\bbQ^t$ and for all
$u\in\bbT$, we have $\bbQ^t\sim_{u,F}\R^u$ and
$\bbQ^t\sim_{u,I}\R^{u^+}$. Then $\bbQ^t\in\Q_1^F\cap\Q^I_2=\Q$ and
hence $\bbQ\in\Q^F$. We make the same argument for the $I$-part.
\end{proof}

\begin{remark}The set $\Q$ defined in Lemma \ref{ca} is given by:
$$
\Q =\left\{\prod\limits_{t\in\bbT}\left(\dfrac{Z^t_{t^+}}{Z^t_{t}}
\times\dfrac{W^{t}_{t+1}}{W^{t}_{t^+}}\right);\;Z^t\in\Q_1,\,W^t\in\Q_2\;\mbox{for
all}\;t\in\bbT\right\}.
$$
\end{remark}

Now we characterize the class, denoted by $\Psi(\Pi)$, of
time-consistent coherent risk measures $\rho$ that satisfy
$\Q^F=\hatpi^F$ with $\Q=\Q^\rho$.

\begin{theorem}
\label{built0}$\rho\in\Psi(\Pi)$ iff there exists some non empty set
$\Phi$ of $\P$-absolutely continuous probabilities measures such
that $\Q=\hatpi^{F}\cap\Phi^{I}$. In this case
$$
\rho_t=\rho^{\hatpi}_t\circ\rho^\Phi_{t^+}\circ\rho_{t+1},
$$
for $t\in\bbT$. In particular if we suppose that $\Pi$ is
time-consistent w.r.t the filtration $(\F_t)_{t\in\bbT}$, then for
all purely financial claims $X\in \Linf(\F)$, we have:
$$
\rho_t(X)=\rho^\Pi_t(X),
$$
for $t\in\bbT$.
\end{theorem}
\begin{proof}Suppose that $\rho\in\Psi(\Pi)$ and define $\Phi=\Q$, we obtain
then by Lemma \ref{bml4kk}, $\Q=\Q^F\cap\Q^I=\Pi^F\cap\Phi^I$.
Conversely suppose that there exists some non empty set
$\Phi\subset\calP$ such that $\Q=\Pi^{F}\cap\Phi^{I}$. From Lemma
\ref{ca} we have $\Q^F=\Pi^{F}$ and $\Q^I=\Phi^{I}$, we deduce that
$\Q$ is time-consistent. To prove the last assertion remark that for
all $t\in\bbT$ and $Y\in\Linf(\F_{t+1})$, we have:
$$
\rho^\Phi_{t^+}(Y)\defto
\mbox{ess-sup}_{\bbQ\in\Phi}\E_\bbQ(Y|\F_{t+1}\otimes\calI_t)=Y,
$$
and
$$
\rho^{\hatpi}_{t}(Y)=\rho^\Pi_{t}(Y).
$$
We prove then by induction on $t=T-1,\ldots,0$ that
$\rho_t(X)=\rho^\Pi_t(X)$ for all $X\in\Linf(\F)$. For $t=T-1$ we
obtain
$$
\rho_{T-1}(X)=
\rho^{\hatpi}_{T-1}\circ\rho^\Phi_{(T-1)^+}(X)=\rho^{\Pi}_{T-1}(X).
$$
Suppose that the induction hypothesis is true until $t+1$, we shall
prove it for $t$. We get
$$
\rho_t(X)=\rho^{\hatpi}_t\circ\rho^\Phi_{t^+}\circ\rho_{t+1}(X)=
\rho^{\hatpi}_t\circ\rho^\Phi_{t^+}\circ\rho^\Pi_{t+1}(X).
$$
Remark that $\rho^\Pi_{t+1}(X)\in\Linf(\F_{t+1})$, then
$$
\rho_t(X)=\rho^{\Pi}_t\circ\rho^\Pi_{t+1}(X)=\rho^{\Pi}_t(X),
$$
from the time-consistency of $\Pi$.
\end{proof}

\end{document}